\documentstyle[amstex,amssymb,12pt,theorem]{article}
\newcommand{\CC}{{\Bbb{C}}}

\newcommand{\RR}{{\Bbb{R}}}
\newcommand{\QQ}{{\Bbb{Q}}}

\newcommand{\ZZ}{{\Bbb{Z}}}
\newcommand{\FF}{{\Bbb{F}}}
\newcommand{\PP}{{\Bbb{P}}}

\theoremheaderfont{\bf\normalshape}\theorembodyfont{\it}
\newtheorem{Thm}{Theorem}
\newtheorem{theorem}{Theorem}[section]
\newtheorem{proposition}[theorem]{Proposition}
\newtheorem{corollary}[theorem]{Corollary}
\newtheorem{lemma}[theorem]{Lemma}
\theorembodyfont{\rm}
 
\newenvironment{Proof}{\begin{ProofwCaption}{Proof}}{\end{ProofwCaption}}
\newenvironment{ProofwCaption}[1]%
  {\addvspace\theorempreskipamount \noindent{\it #1.}\rm}%
  {\qed \par \addvspace\theorempostskipamount}
\newcommand{\qedsymbol}{\mbox{$\Box$}}
\newcommand{\qed}{\quad\qedsymbol}
\begin{document}
\title{The Fundamental Group of Some Siegel \\
Modular Threefolds}
\author{Klaus Hulek and G. K. Sankaran}
\date{}
\maketitle
The moduli space ${\cal A}_{1,p}$ of abelian surfaces with a polarisation
of type
$(1,p)$ and a level structure, for $p$ an odd prime, is a singular
quasi-projective variety. It has been studied in detail in [HKW1]. In the
present note we prove the following.
\begin{Thm}\label{0:t 0.1}
Let X be any desingularisation of an algebraic compactification of ${\cal
A}_{1,p}$
for $p$ an odd prime. Then $X$ is simply connected.
\end{Thm}
In fact it is enough to prove this for one such $X$, since any two $X$
have the same fundamental group.\\
If $p=5$ or $7$ then $X$ is known to be rational and therefore simply
connected. For $p=5$ this goes back to [HM]; for $p=7$ it has recently
been shown in [MS] that $X$ is birationally equivalent to a Fano variety
of type $V_{22}$, which is known to be rational. It does not seem to be
known whether $X$ is rational if $p=3$, but
it is shown in [HS] that if $p$ is large then $X$ is of general type.

\section{Generalities about fundamental groups}
The general facts in this section are all well known.
\begin{lemma}\label{1:l 1.1}
Let $M$ be a connected simply connected real manifold and $G$ a group acting
discontinuously on $M$. Let $x\in M$ be a base point. Then the quotient
map $\varphi:M\longrightarrow M/G$ induces a map $\psi:G\longrightarrow
\pi_1(M/G,\varphi(x))$ which is surjective.
\end{lemma}
Note that we do not require the action of $G$ to be free or even
effective. The map $\psi$ is defined as follows. If $g\in G$ let
$\theta_g:[0,1]\longrightarrow M$ be any path in $M$ from $x$ to $g(x):$
that is, $\theta_g(0)=x$ and $\theta_g(1)=g(x)$. Then
$\varphi\circ\theta_g:[0,1]\longrightarrow M/G$ is a closed loop in $M/G$,
with $\varphi\circ\theta_g(0)=\varphi\circ\theta_g(1)=\varphi(x)$. We put
$\psi(g)=[\varphi\circ\theta_g] \in \pi_1(M/G,\varphi(x))$. One checks easily
that $\psi$ is well defined and that it is a group homomorphism. It is
shown in [G] that $\psi$ is surjective.
\begin{lemma}\label{1:l 1.2}
Let $M$ be a connected complex manifold and $M_0$ an analytic subvariety.
Then the inclusion
induces a map $\pi_1(M\setminus M_0,x)\longrightarrow \pi_1(M,x)$ for any base
point
$x\in M$ which is surjective for $\operatorname{codim }_{ \CC} M_0\ge 1$ and an
isomorphism for $\operatorname{codim }_{ \CC} M_0 \ge 2$.
\end{lemma}
\begin{lemma}\label{1:l 1.3}
If $X_1$ and $X_2$ are birationally equivalent smooth projective varieties then
$\pi_1(X_1)\cong \pi_1(X_2)$.
\end{lemma}
Proofs of both these lemmas can be found in [HK].

\section{Applications to the case of ${\cal A}_{1,p}$}
We fix some notation. $S_2$ is the Siegel upper half-plane of degree 2. We
denote by
${\cal A}_{1,p}$ the moduli space of abelian surfaces with a
$(1,p)$-polarisation and a level
structure, and by ${\cal A}(p^2)$ the moduli space of abelian surfaces of
level $p^2$. We let
$\Gamma_{1,p}$ and $\Gamma(p^2)$ be the corresponding arithmetic subgroups of
$\operatorname{Sp}(4,\ZZ)$, so
$$
\Gamma_{1,p}=\left\{\gamma\in
\operatorname{Sp}(4,\ZZ)|\gamma-{\mbox{\boldmath$1$}}\in\left(
\begin{array}{cccc}
\ZZ   &   \ZZ   &   \ZZ   & p\ZZ\\
p\ZZ  &   p\ZZ  &   p\ZZ  & p^2\ZZ\\
\ZZ   &   \ZZ   &   \ZZ   & p\ZZ\\
\ZZ   &   \ZZ   &   \ZZ   & p\ZZ\\
\end{array}
\right) \right\}
$$
(see [HKW1]) and $\Gamma (p^2)$ is the principal congruence subgroup of
level $p^2$. We
denote by $\Gamma^0_{1,p}$ the arithmetic subgroup of
$\operatorname{Sp}(4,\QQ)$

$$
\Gamma_{1,p}^0=\left\{\gamma\in \operatorname{Sp}(4,\QQ)|\gamma\in\left(
\begin{array}{cccc}
\ZZ   &   \ZZ   &   \ZZ   & p\ZZ\\
p\ZZ  &   \ZZ   &   p\ZZ  & p\ZZ\\
\ZZ   &   \ZZ   &   \ZZ   & p\ZZ\\
\ZZ   &   \frac 1 p\ZZ   &   \ZZ   & \ZZ\\
\end{array}
\right) \right\}
$$
(in terms of moduli, this corresponds to keeping the $(1,p)$-polarisation
but forgetting
the level structure). We denote by ${\cal A}^*_{1,p}$ the toroidal
compactification of
${\cal A}_{1,p}$ described in [HKW1], and by ${\cal A}^*(p^2)$  the usual
toroidal (or Igusa)
compactification of ${\cal A}(p^2)$. Since $\Gamma(p^2)$ is neat, ${\cal
A}^*(p^2)$ is
nonsingular: indeed, any nonsingular compactification of ${\cal A}(p^2)$
will serve our purpose.\\
Let ${\cal A}'_{1,p}$ be the smooth part of ${\cal A}_{1,p}$. The space
${\cal A}_{1,p}$ is a
quotient of $S_2$ by $\Gamma_{1,p}$ and ${\cal A}'_{1,p}$ is the quotient by
$\Gamma_{1,p}$ of an open part $S'_2 \subset S_2$. Since the singular locus
of ${\cal
A}_{1,p}$ has codimension $2,S_2'$ is the complement of an analytic subset
of codimension
$2$ in $S_2$. In particular, by Lemma 1.2, $S'_2$ is simply connected.\\
{\em Remark. }
The quotient map $S'_2 \longrightarrow {\cal A}'_{1,p}=S'_2/\Gamma_{1,p}$ is
not
unramified. The ramification is in codimension $1$ but the isotropy groups
are generated by
reflections so the quotient is smooth. If we restrict to the complement of
the ramification
locus our covering space will not necessarily be simply connected.
\begin{proposition}\label{2:p 2.2}
Let $f:X \longrightarrow {\cal A}^*_{1,p}$ be a resolution of singularities
(so $X$ is smooth
and $f$ is a birational morphism which is an isomorphism away from
$f^{-1}(\operatorname{Sing} {\cal A}^*_{1,p}))$. Then there are surjective
homomorphisms
$\Gamma_{1,p} \longrightarrow \pi_1({\cal A}'_{1,p})$ and $\pi_1({\cal
A}'_{1,p})
\longrightarrow \pi_1(X)$. \end{proposition}
\begin{Proof}
The map $\Gamma_{1,p} \longrightarrow \pi_1({\cal A}'_{1,p})$ exists and is
surjective by
Lemma 1.1. The map $\pi_1({\cal A}'_{1,p})\longrightarrow \pi_1(X)$ comes
from Lemma 1.2:
we may consider ${\cal A}'_{1,p}$ as a subset of $X$ via $f^{-1}$, which
exists on ${\cal
A}'_{1,p}$.
\end {Proof}
We now choose a convenient resolution of singularities to work with.
$\Gamma_{1,p}$ is a
normal subgroup of $\Gamma^0_{1,p}$ and the quotient (which is isomorphic to
$\mbox{SL}(2,\FF_p)$ as an abstract group) acts on ${\cal A}^*_{1,p}$ (this
is shown in
[HKW1]). By the general results of Hironaka ([H]) or by an easy explicit
construction we can
choose $f:X\longrightarrow {\cal A}^*_{1,p}$ to be an equivariant
resolution, so that
$\Gamma^0_{1,p}/\Gamma_{1,p}$ acts on $X$. Henceforth $X$ will always be such a
resolution.\\ Let $\psi:\Gamma_{1,p}\longrightarrow \pi_1(X)$ be the
composite of the two
maps in Proposition 2.1. Then $\psi$ is surjective. We shall show that
$\operatorname{Ker}
 \psi=\Gamma_{1,p}$ and hence $\pi_1(X)=1$.

\section{Level $p^2$}
In this section we use results of Kn\"oller ([K]) to prove that
$\Gamma(p^2)$, which is a
subgroup of $\Gamma_{1,p}$, lies in the kernel of $\psi$.\\
Since $\Gamma(p^2)$ is a normal subgroup of $\Gamma_{1,p}$, there is an
action of
$\Gamma_{1,p}/\Gamma (p^2)$ on ${\cal A}(p^2)$ and the quotient is ${\cal
A}_{1,p}$. The
quotient map induces a rational map $h:{\cal A}^*(p^2)--\rightarrow X$,
which is a
morphism over ${\cal A}'_{1,p}$. We can resolve the singularities of this
map so as to get a
diagram
$$
\unitlength1cm
\begin {picture}(5,3)
\put(0.5,2){\vector(0,-1){1}}
\put(0.8,2){\vector(2,-1){2.1}}
\multiput(1.5,0.5)(0.428,0){4}{\line(1,0){0.214}}
\put(2.78,0.5){\vector(1,0){0.214}}
\put(0.5,2.2){$Y$}
\put(0.1,0.4){${\cal A}^*(p^2)$}
\put(3.1,0.4){$X$}
\put(0.2,1.5){$\sigma$}
\put(2.1,0.7){$h$}
\put(2.2,1.6){$\tilde{h}$}
\end{picture}
$$
where $Y$ is smooth, $\sigma$ is an isomorphism over ${\cal
A}'(p^2)=h^{-1}({\cal A}'_{1,p})$
and $\tilde h$ is a morphism. In particular we can think of ${\cal
A}'(p^2)$ as a subset of $Y$
and ${\cal A}'_{1,p}$ as a subset of $X$, and then $\tilde h|_{{\cal
A}'(p^2)}$ is the quotient
map ${\cal A}'(p^2)\longrightarrow {\cal A}'_{1,p}$.
\begin{theorem}\label{3:t 3.1} {\em (Kn\"oller [K]) }
Y  is simply connected.
\end{theorem}
\begin{corollary}\label{3:c 3.2}
$\Gamma(p^2)$ is contained in $\operatorname{Ker} \psi$.
\end{corollary}
\begin{Proof}
The quotient map $q:S'_2\longrightarrow {\cal A}'(p^2)$ induces a
surjection $\Gamma
(p^2)\longrightarrow \pi_1(Y)$. Let $M \in \Gamma(p^2)$, let $x\in S'_2$ be
a base point
and let $\theta_M$ be a path from $x$ to $M(x)$ in $S'_2$. Then there is a
null homotopy
$H_M:[0,1]\times[0,1]\longrightarrow Y$ such that $H_M(0,t)=q \theta_M(t)$ and
$H_M(1,t)=q(x)\in Y$. So $\tilde{h} H_M$ is a null homotopy of $\tilde{h}q
\theta_M$ in
$X$. But $[\tilde{h}q \theta_M]$ is the class $\psi(M)\in \pi_1(X)$, since
$\tilde{h}
q$ is the same as the quotient map $S'_2\longrightarrow {\cal A}'_{1,p}$.
\end {Proof}

\section{Another element of the kernel}
\begin{proposition}\label {4:l 1.1}
The element
$$
M_0=\left (
\begin{array}{cccc}
1  &  0  &  1  &  0\\
0  &  1  &  0  &  0\\
0  &  0  &  1  &  0\\
0  &  0  &  0  &  1
\end{array}
\right )
$$
of $\Gamma_{1,p}$ is in $\operatorname{Ker }\psi$.
\end{proposition}
\begin{Proof}
We follow the procedure of [K], and choose a path
$\theta_{M_0}:[0,1]\longrightarrow S_2$,
namely $$
\theta_{M_0}(t)=ic (1-t){\mbox{\boldmath$1$}} + ict M_0({\mbox{\boldmath$1$}})
$$
where $c\in \RR$ is a positive constant to be specified later.
$\theta_{M_0}$ gives rise to a
closed loop in ${\cal A}_{1,p}$, starting and finishing at
$ic{\mbox{\boldmath$1 $} }$ (so when
we specify $c$ we are just choosing a base point).\\
The map $e_{l_0}:S_2\longrightarrow\CC^*\times \CC \times S$, corresponding
to the central
boundary component $D_{l_0}$ of ${\cal A}^*_{1,p}$ (see [HKW1]) is given by
$$
e_{l_0}:\left (\begin{array}{cc}
\tau_1   &   \tau_2\\
\tau_2   &   \tau_3
\end{array}
\right ) \longmapsto (e^{2\pi i \tau_1}, \tau_2,\tau_3).
$$
So $ e_{l_0}\theta_{M_0}(t)=(e^{-2\pi c} e^{-2\pi i t}, 0, ic)$. The homotopy
$$
H(s,t)=((se^{2\pi i t}+1-s)e^{-2\pi c}, 0, ic)
$$
is evidently a null homotopy of $e_{l_0}\theta_{M_0}$ in $\CC \times \CC
\times S_1$. As in [K], we
let $V$ be the interior of the closure of $S_2/P_{l_0}$ in $\CC \times \CC
\times S_1:$
then the image of $H$ lies in $V$. But by [HKW2], ${\cal A}^*_{1,p}$ is
nonsingular near
$D_{l_0}:$ the map $V\longrightarrow {\cal A}^*_{1,p}$ is branched but in
some open set $U
\subset V$ the branching comes from reflections. If we choose $c$
sufficiently large then the
image of $H$ will lie entirely in $U$, and $H$ will therefore give rise to
a null homotopy
which does not meet  the singular locus of ${\cal A}^*_{1,p}$. In
particular, the loop in ${\cal
A}_{1,p}$ corresponding to $\theta_{M_0}$ actually lies in ${\cal
A}'_{1,p}$, and $\theta_{M_0}$
is actually a path in $S'_2$. So $\psi(M_0)$ is trivial (it is even trivial as
an element of
$\pi_1({\cal A}^*_{1,p}\setminus \mbox{Sing }{\cal A}^*_{1,p}))$.
\end{Proof}
Now we know that $\operatorname{Ker} \psi$ contains $M_0$ and all of
$\Gamma(p^2)$. It is also, of
course, a normal subgroup of $\Gamma_{1,p}:$ but because we could choose
$X$ to have an
action of $\Gamma^0_{1,p}/\Gamma_{1,p}$ it is actually normal as a subgroup of
$\Gamma^0_{1,p}$, since $\Gamma^0_{1,p}$ acts on
$\pi_1(X)=\Gamma_{1,p}/\operatorname{Ker}\psi$ by
conjugation in $\mbox{Sp}(4,\QQ)$.
\section{Calculations in $\Gamma _{1,p}$}
The remark at the end of section 4, above, shows that to prove the theorem
it is enough to
check that any normal subgroup of $\Gamma^0_{1,p}$ containing $\Gamma(p^2)$
and $M_0$
must contain $\Gamma_{1,p}$. We shall give a set of generators for
$\Gamma_{1,p}$ and show
how to produce each one, starting with $M_0$ and $\Gamma(p^2)$ and using
multiplication,
inversion, and conjugation by elements of $\Gamma^0_{1,p}$. A similar
procedure is
carried out by Mennicke [M] for principal congruence subgroups of
$\mbox{Sp}(2n,\ZZ)$.
\begin{proposition}\label{5:p1.11}
$\Gamma_{1,p}$ is generated by the elements $M_1, M_2, M_3$ and $M_4$
together with the
subgroups $j_1 (\operatorname{SL}(2,\ZZ))$ and $j_2(\Gamma_1(p))$, where
$$
M_1=\left (
\begin{array}{cccc}
1  &  0  &  0  &  0\\
0  &  1  &  0  &  0\\
0  &  1  &  1  &  0\\
1  &  0  &  0  &  0\\
\end{array}\right ) \qquad
M_2=\left (
\begin{array}{cccc}
1  &  0  &  0  &  p\\
0  &  1  &  p  &  0\\
0  &  0  &  1  &  0\\
0  &  0  &  0  &  1\\
\end{array}\right )
$$
$$
M_3=\left (
\begin{array}{cccc}
1  &  1  &  0  &  0\\
0  &  1  &  0  &  0\\
0  &  0  &  1  &  0\\
0  &  0  & -1  &  1
\end{array}\right )\qquad
M_4=\left (
\begin{array}{cccc}
1  &  0  &  0  &  0\\
-p &  1  &  0  &  0\\
0  &  0  &  1  &  p\\
0  &  0  &  0  &  1
\end{array}\right )
$$
$$
j_1(\operatorname{SL}(2,\ZZ))=\left\{ \left (
\begin{array}{cccc}
a  &  0  &  b  &  0\\
0  &  1  &  0  &  0\\
c  &  0  &  d  &  0\\
0  &  0  &  0  &  1\\
\end{array}
\right )
|\left (
\begin{array}{cc}
a  &  b\\
c  &  d
\end{array}
\right ) \in \operatorname{SL}(2,\ZZ) \right \}
$$
$$
j_2(\Gamma_1(p))=\left\{ \left (
\begin{array}{cccc}
1  &  0  &  0  &  0\\
0  &  a  &  0  &  pb\\
0  &  0  &  1  &  0\\
0  & c/p &  0  &  d
\end{array}
\right )
|\left(
\begin{array}{cc}
a  &  b\\
c  &  d
\end{array}
\right ) \in
\Gamma_1(p) \right \}.
$$
\end{proposition}
{\em Remark. }
Of course one could find a finite set of generations for $\Gamma_{1,p}$ by
doing so for
$\mbox{SL}(2,\ZZ)$ and for $\Gamma_1(p)$.

\begin{Proof}
It is easier to work, as in Chapter 1 of [HKW1], with
$$
\tilde{\Gamma}_{1,p}=\left \{ \gamma \in \mbox{Sp} (\Lambda, \ZZ)| \gamma
\equiv \left (
\begin{array}{cccc}
\ast  &  \ast  &  \ast  &  \ast\\
0      &  1       &  0       &  0\\
\ast  &  \ast  &  \ast  &  \ast\\
0      &  0       &  0       &  1
\end{array}
\right ) \mbox{mod } p \right \}
$$
where $\Lambda$ is the symplectic form
$$ \left (
\begin{array}{cc}
\phantom{-}0   &   \begin{array}{cc}
1 & 0\\
0 & p
\end{array}\\
\begin{array}{cc}
-1                   & \phantom{-}0\\
\phantom{-}0  & -p
\end{array} & 0
\end{array}
\right ).
$$
$\tilde{\Gamma}_{1,p}$ is an $\mbox{Sp}(4,\QQ)$-conjugate of
$\Gamma_{1,p}$. More
precisely $$
\tilde{\Gamma}_{1,p}= R {\Gamma}_{1,p} R^{-1}
$$
where
$$
R=\left (
\begin{array}{cccc}
1 &&&\\
& 1 &&\\
&& 1 &\\
&&& p
\end{array}
\right ).
$$
Thus we wish to show that $\tilde{\Gamma}_{1,p}$ is generated by
$$
\tilde{M}_1=\left (
\begin{array}{cccc}
1  &  0  &  0  &  0\\
0  &  1  &  0  &  0\\
0  &  1  &  1  &  0\\
p  &  0  &  0  &  1\\
\end{array}\right ),\qquad
\tilde{M}_2=\left (
\begin{array}{cccc}
1  &  0  &  0  &  1\\
0  &  1  &  p  &  0\\
0  &  0  &  1  &  0\\
0  &  0  &  0  &  1\\
\end{array}\right ) ,
$$
$$
\tilde{M}_3=\left (
\begin{array}{cccc}
1  &  1  &  0  &  0\\
0  &  1  &  0  &  0\\
0  &  0  &  1  &  0\\
0  &  0  &  -p &  1\\
\end{array}\right ),\qquad
\tilde{M}_4=\left (
\begin{array}{cccc}
1  &  0  &  0  &  0\\
-p &  1  &  0  &  0\\
0  &  0  &  1  &  1\\
0  &  0  &  0  &  1\\
\end{array}\right ),
$$
$$
\tilde{\jmath}_1(\mbox{SL}(2,\ZZ))=\left \{
\left (
\begin{array}{cccc}
a  &  0  &  b  &  0\\
0  &  1  &  0  &  0\\
c  &  0  &  d  &  0\\
0  &  0  &  0  &  1\\
\end{array}\right ) | \left (
\begin{array}{cc}
a & b\\
c & d
\end{array} \right ) \in \mbox{SL}(2,\ZZ)\right \}
$$
and
$$
\tilde{\jmath}_2(\Gamma_1(p))=\left \{
\left (
\begin{array}{cccc}
1  &  0  &  0  &  0\\
0  &  a  &  0  &  b\\
0  &  0  &  1  &  0\\
0  &  c  &  0  &  d\\
\end{array}\right ) |
\left (
\begin{array}{cc}
a & b\\
c & d
\end{array} \right ) \in \Gamma_1(p) \right \}.
$$
Recall that a vector $v \in \ZZ^4$ is short if and only if there exists $w
\in \ZZ^4$ such
that $v \Lambda ^t w =1:$ otherwise it is long. If $K \in
\tilde{\Gamma}_{1,p}$ then the
second row of $K$ is certainly long. Suppose the first row is long also.
Then, since the
first two rows of K span a $\Lambda$-isotropic plane, $K_{14}\equiv 0
\mbox{ mod } p$. But
then $$
K\equiv\left (
\begin{array}{cccc}
0       &  \ast  &  0      &  0\\
0       &  1      &  0       &  0\\
\ast  &  \ast  &  \ast  &  \ast\\
\ast  &  \ast  &  \ast  &  \ast\\
\end{array}\right ) \mbox{mod } p
$$
so $K$ is not invertible mod $p$. But $K \in \mbox{Sp}(4,\ZZ)$ and in
particular $K$ is
invertible.\\
So the first row of $K$ must be short. It is shown in the course of the proof
of
Proposition 3.38 in [HKW1] that, by multiplying $K$ on the right by a
suitable product of
$\tilde{M}_i$ and elements of $\tilde{\jmath}_1(\mbox{SL}(2,\ZZ))$, we may
assume that
the first row of $K$ is $(1, 0, 0, 0)$. Then the symplectic condition gives
$$
K=\left (
\begin{array}{cccc}
1       &  0      &  0      &  0\\
\ast   &  a      &  0       &  c\\
\ast  &  \ast  &  1  &  \ast\\
\ast  &   b      &  0  &  d
\end{array}\right )
$$
for some
$\left(
\begin{array}{cc}
a & b\\
c & d
\end{array}\right) \in \Gamma_1 (p)$. Multiplying by $\tilde {\jmath}_2\left(
\begin{array}{cc}
a & b\\
c & d
\end{array}\right)^{-1}$ transforms this into
$$
\left (
\begin{array}{cccc}
1       &  0      &  0   &  0\\
-np   &   1      &  0   & 0\\
\ast  &   m     &  1  &  n\\
mp    &   0      &  0  &  1
\end{array}\right )
$$
for some integers $m,n$. Multiplying on the right by $\tilde{M}_4^{-n}
\tilde{M}_1^{-m}$
reduces to
$$
\left (
\begin{array}{cccc}
1       &  0      &  0   &  0\\
0      &   1      &  0   &  0\\
\ast  &   0     &  1   &  0\\
0      &   0      &  0  &  1
\end{array}\right )
$$
which is in $\tilde{\jmath}_1(\mbox{\mbox{SL}}(2,\ZZ))$. This proves the
proposition.
\end{Proof}
Now we produce all the generators using the allowable means described above.
\begin{proposition}\label{5:p 5.2}
$M_1, M_2, M_3, M_4, j_1({\operatorname{SL}}(2,\ZZ))$ and
$j_2(\Gamma_1(p))$ can all
be generated from $M_0$ and $\Gamma(p^2)$ by multiplication, inversion, and
conjugation in
$\Gamma^0_{1,p}$.
\end{proposition}
\begin{Proof}
(i)\quad $M_0=j_1 \left( \begin{array}{cc}
1  &  1\\
0  &  1
\end{array} \right )  $
and the smallest normal subgroup of $\mbox{SL}(2,\ZZ)$ containing $M_0$ is
the whole of
$\operatorname{SL}(2,\ZZ)$, e.g. [B]. Hence we can make
$j_1(\operatorname{SL}(2,\ZZ))$.\\
(ii)\quad $ M^{-1}_4 M_0 M_4 M_0^{-1}= \left (
\begin{array}{ccc|c}
& {\mbox{\boldmath$1$}} & &
\begin{array}{cc}
0  &  p\\
p  &  p^2
\end{array} \\
\hline
& 0 \vphantom{\begin{array}{cc}
0  &  p\\
p  &  p^2
\end{array} } &  & {\mbox{\boldmath$1$}}
\end{array} \right ) $
and $L_1=
 \left (
\begin{array}{ccc|c}
 &{\mbox{\boldmath$1$}}&  &
\begin{array}{cc}
0  &  0\\
0  &  p^2
\end{array} \\
\hline
&0\vphantom{\begin{array}{cc}
0  &  p\\
p  &  p^2
\end{array} } &   & {\mbox{\boldmath$1$}}
\end{array} \right )
\in \Gamma(p^2)$, so \\
$M_2=M_4^{-1} M_0 M_4 M_0^{-1} L_1^{-1}$ and we have made $M_2$.\\
(iii)\quad $ M_1 M_0 M_1^{-1} M_1^{-1} M_0 M_1=
\left (
\begin{array}{cccc}
1       &  0      &  2   &  0\\
0      &   1      &  0   &  0\\
0      &   0     &  1   &  0\\
0      &  -2      &  0  &  1
\end{array}\right )
$
and since we got $j_1(\operatorname{SL}(2,\ZZ))$ in (i) we may multiply by
$j_1\left (
\begin{array}{cc} 1  &  -2\\
0  &   1
\end{array}\right )$ to get
$$
\left (
\begin{array}{cccc}
1       &  0      &  0   &  0\\
0      &   1      &  0   &  0\\
0      &   0      &  1   &  0\\
0      &  -2      &  0  &  1
\end{array}\right )
$$
Call this matrix $L_2$. Since $p$ is odd we can find integers $\lambda,
\mu$  such that
$-2\lambda + p^2 \mu=1$. Consider the matrix  $L_3= j_2\left (
\begin{array}{cc}
1  &  0\\
p^3  &   1
\end{array}
\right ) $,
so $L_3\in \Gamma(p^2)$. We then have
$$
L^{\lambda}_2 L^{\mu}_3 = \left (
\begin{array}{cccc}
1       &  0      &  0   &  0\\
0      &   1      &  0   &  0\\
0      &   0     &   1   &  0\\
0      &   1     &   0   &  1
\end{array}\right ) = j_2  \left (\begin{array}{cc}
1  &  0\\
p  &   1
\end{array}
\right).
$$
Call this element $L_4$. We also have $j_2 \left (\begin{array}{cc}
1  &  p\\
0  &  1
\end{array}
\right) \in \Gamma (p^2)$, but we have not got all of $j_2(\Gamma_1(p))$ yet.
However we need $L_4$ as an auxiliary.\\
(iv)\quad
$ L_4 M_1 M_0 M_1^{-1} M_0^{-1}=M_3^{-1}$ so we can now make $M_3.$\\
(v)\quad
Put $L_5=j_1\left (\begin{array}{cc}
1 & 0\\
1 & 1
\end{array} \right ).$ Then
$$
L_5^{-1} M_2^{-1} L_5 M_2 L_1^{-1} = M_4
$$
and we have made $M_4$.\\
(vi)\quad
$M_3 L_5 L_4 M_3^{-1} L_5^{-1} L_2=M_1^{-1}$, so we can make $M_1$.\\
(vii)\quad
It remains to show how to get all of $j_2(\Gamma_1(p))$. So far we have not
used the freedom
to conjugate by $\Gamma^0_{1,p}$ rather than just by $\Gamma _{1,p}$. Define
$$
\Gamma_1'(p^2)=\left\{
Q\in\Gamma_1(p)|Q-{\mbox{\boldmath$1$}}=\left(\begin{array}{cc}
p^2\ZZ   &   p\ZZ\\
p^3\ZZ   &   p^2\ZZ
\end{array}\right) \right\}.
$$
Thus $Q\in\Gamma'_1(p^2)$ if and only if $j_2(Q)\in
\Gamma(p^2)$. If instead $Q\in \mbox{SL}(2,\ZZ)$ then $j_2(Q)\in
\Gamma^0_{1,p}$,
and we have already generated $j_2(P)$, where $P=\left(\begin{array}{cc}
1  &  0\\
p  &  1
\end{array}\right )$, so the problem is now to generate $\Gamma_1(p)$ using
$P$, elements
of $\Gamma'_1(p^2)$, and conjugation by elements of $\mbox{SL}(2,\ZZ)$.\\
A general element of $\Gamma _1(p)$ is of the form $\left( \begin{array}{cc}
\lambda p+1   &   \alpha p\\
\beta p   &   \mu p+1
\end{array}\right )$ and since it has determinant equal to 1
$$
(\lambda + \mu)p+(\lambda\mu - \alpha\beta)p^2=0,
$$
and in particular $\lambda + \mu\equiv 0 \mbox{ mod } p$. Suppose that
$\lambda\equiv
0\mbox{ mod } p$, so that $\mu\equiv 0\mbox{ mod } p$ also. Then
$$
P^{-\beta} \left ( \begin{array}{cc}
\lambda p+1   &   \alpha p\\
\beta p           &   \mu p+1
\end{array}\right ) = \left (\begin{array}{cc}
\lambda p+1            &   \alpha p\\
-\beta\lambda p^2   &   -\alpha\beta p^2+\mu p+1
\end{array}\right )
$$
which is in $\Gamma_1'(p^2)$. So we can generate any element of
$\Gamma_1(p)$ for which
$\lambda\equiv 0\mbox{ mod } p$.\\
Suppose then that $(\lambda, p)=1$. If $(\alpha, p)=1$ also then
$p|(\lambda-k\alpha)$ for
some integer $k$. But now we can conjugate by $\left ( \begin{array}{cc}
1  &  0\\
k  &  1
\end{array}\right ) \in \mbox{SL}(2,\ZZ)$:
$$
\left( \begin{array}{cc}
1   &   0\\
k   &   1
\end{array}\right)
\left ( \begin{array}{cc}
\lambda p+1   &   \alpha p\\
\beta p           &   \mu p+1
\end{array}\right )
\left( \begin{array}{cc}
1   &   0\\
-k   &   1
\end{array}\right)=
\left ( \begin{array}{cc}
(\lambda -k\alpha)p+1   &   \alpha p\\
\beta' p           &   (\mu+k\alpha) p+1
\end{array}\right )
$$
and so we can generate all elements of $\Gamma_1(p)$ for which
$\lambda\equiv0\mbox{ mod } p$ or $ \alpha\not\equiv 0 \mbox{ mod } p$. The
remaining
elements are of the form
$$
\left( \begin{array}{cc}
\lambda p+1   &   \alpha' p^2\\
\ast   &   \mu p+1
\end{array}
\right )
$$
and we multiply by $\left( \begin{array}{cc}
1   &   p\\
0   &   1
\end{array}\right) \in \Gamma'_1(p^2)$:
$$
\left( \begin{array}{cc}
1   &   p\\
0   &   1
\end{array}\right)
\left( \begin{array}{cc}
\lambda p+1   &   \alpha' p^2\\
\ast   &   \mu p+1
\end{array}
\right )=
\left( \begin{array}{cc}
\ast   &   (\alpha'+\mu) p^2+p\\
\ast   &   \mu p+1
\end{array}
\right )
$$
and $(\alpha'+\mu)p+1\not\equiv 0 \mbox { mod } p$ so this is something we have
already
generated.
\end{Proof}

\section*{References}
\begin{enumerate}
 \item [{[B]}] J.L. Brenner, The linear homogeneous group III. Ann. Math
{\bf71} (1960),
210-223.
\item[{[G]}] J. Grosche, \"Uber die Fundamentalgruppen von Quotientr\"aumen
Siegelscher Modulgruppen, J. reine angew. Math. {\bf 281} (1976), 53-79.
\item[{[HK]}] H. Heidrich \& F.W. Kn\"oller, \"Uber die Fundamentalgruppen
Siegelscher
Modulvariet\"aten vom Grade 2, Manuscr. Math. {\bf 57} (1987), 249-262.
\item[{[H]}] H. Hironaka, Resolution of singularities of an algebraic
variety over a field of
characteristic zero, I, II, Ann. Math. {\bf 79} (1964), 109-326.
\item[{[HM]}] G. Horrocks \& D. Mumford, A rank 2 vector bundle on $\PP^4$
with 15,000
symmetries, Topology {\bf 12} (1973), 63-81.
\item[{[HKW1]}] K. Hulek, C. Kahn \& S.H. Weintraub, Moduli spaces of
abelian surfaces:
compactification, degenerations and theta functions, de Gruyter, Berlin 1993.
\item[{[HKW2]}] K. Hulek, C. Kahn \& S.H. Weintraub, Singularities of the
moduli space of
certain Abelian surfaces, Comp. Math. {\bf 79} (1991), 231-253.
\item[{[HS]}] K. Hulek, G.K. Sankaran, The Kodaira dimension of certain
moduli spaces of
Abelian surfaces, Comp. Math., to appear.
\item[{[K]}] F.W. Kn\"oller, Die Fundamentalgruppen der Siegelschen
Modulvariet\"aten, Abh.
Math. Sem. Univ. Hamburg {\bf 57} (1987), 203-213.
\item[{[MS]}] N. Manolache \& F. Schreyer, preprint 1993.
\item[{[M]}] J. Mennicke, Zur Theorie der Siegelschen Modulgruppe, Math.
Annalen {\bf 159}
(1965), 115-129.
\end{enumerate}
\begin{tabular*}{13cm}{l@{\extracolsep{\fill}}l}
K. Hulek                            & G.K. Sankaran\\
Institut f\"ur Mathematik & Department of Pure Mathematics\\
Universit\"at Hannover     & and Mathematical Statistics\\
Postfach 6009                  & University of Cambridge\\
D-30060 Hannover            & Cambridge CB2 1SB\\
Germany                           & England
\end{tabular*}
\end{document}